\documentclass[twoside,11pt]{article}

\usepackage{amssymb}

\renewcommand{\k}{\Bbbk}

\newcommand{\N}{{\mathbb N}}
\newcommand{\Z}{{\mathbb Z}}
\newcommand{\Q}{{\mathbb Q}}

\newcommand{\M}{\mathcal M}

\newcommand{\I}{\mathcal I}

\newcommand{\GKdim}{{\rm GKdim}}

\newcommand{\im}{{\rm im}}

\newcommand{\rk}{{\rm rk}}
\newcommand{\st}{{\rm st}}

\newcommand{\tq}{\,|\,}
\newcommand{\qed}{\hfill \rule{1.5mm}{1.5mm}}
\newcommand{\proof}{\noindent {\it Proof: $\;$}}

\newtheorem{theorem}{Theorem}[section]
\newtheorem{proposition}[theorem]{Proposition}

\newtheorem{definition}[theorem]{Definition}
\newtheorem{lemma}[theorem]{Lemma}

\newtheorem{remark}[theorem]{Remark}

\newtheorem{subtheorem}{Theorem}[subsection]
\newtheorem{subproposition}[subtheorem]{Proposition}

\newtheorem{subdefinition}[subtheorem]{Definition}
\newtheorem{sublemma}[subtheorem]{Lemma}
\newtheorem{subexample}[subtheorem]{Example}

\newtheorem{subcorollary}[subtheorem]{Corollary}
\newtheorem{subremark}[subtheorem]{Remark}

\newcommand{\titre}{Quantum analogues of Schubert varieties\\ 
in the grassmannian.}

\begin{document}

\title{{\vspace{-1.5cm} \bf \titre}}
\author{T.H. Lenagan and L. Rigal\footnote{This research 
was supported by the Leverhulme Research 
Interchange Grant
F/00158/X.}}
\date{}
\maketitle

\begin{abstract} 
We study quantum Schubert varieties from the point of view of regularity
conditions. More precisely, we show that these rings are 
domains which are maximal orders and are AS-Cohen-Macaulay and we determine which of them are
AS-Gorenstein. One key fact that enables us to prove these results is that
quantum Schubert varieties are quantum graded algebras with a straightening
law that have a unique minimal element in the defining poset. We prove a
general result showing when such quantum graded algebras are maximal orders. Finally, we exploit these results to show that
quantum determinantal rings are maximal orders.
\end{abstract}



\section*{Introduction.}

Since 
the appearance of quantum groups in the eighties, there have been several
attempts to define quantum analogues of coordinate rings of grassmannian
varieties and, more generally, of flag varieties. Here, we are interested in
such deformations for grassmannian varieties and we follow the approach of
Lakshmibai and Reshetikhin (see [LakRe]). Hence, we start with the (usual)
quantum deformation, denoted by 
${\mathcal O}_q(M_{m,n}(\k))$, of the coordinate ring of $m
\times n$ matrices. Then, denoting by $G_{m,n}(\k)$ the grassmanian of
$m$-dimensional subspaces in $\k^n$, the quantum deformation of the coordinate
ring of $G_{m,n}(\k)$ that we consider is the $\k$-subalgebra of
${\mathcal O}_q(M_{m,n}(\k))$ generated by the maximal quantum minors. We denote 
this algebra by
${\mathcal O}_q(G_{m,n}(\k))$ and call it the {\em quantum grassmannian} for simplicity.
Precise definitions are recalled in Section \ref{sec-basic-definitions}.

Our main interest, in this paper, is the study of a family of quotients of
${\mathcal O}_q(G_{m,n}(\k))$ which appear in [LakRe] and are natural quantum analogues
of coordinate rings of Schubert varieties in $G_{m,n}(\k)$. These {\em quantum
Schubert varieties} have already been studied, to some extent, in [LenRi(2)].
There, they were used as a central tool to show that ${\mathcal O}_q(G_{m,n}(\k))$ is a
quantum graded algebra with a straightening law. Details on the notion of
quantum graded algebra with a straightening law (quantum graded A.S.L. for
short) can be found in Section \ref{q-gr-ASL} below and in [LenRi(2)] where 
this notion 
was introduced. Roughly speaking, such an algebra is endowed with a standard
monomial basis build on the elements of a finite 
partially ordered subset of generators.
In addition, the quantum graded A.S.L. structure provides a good control on
the way such standard monomials multiply with each other by means of the
so-called {\em straightening law} and {\em commutation law}.

In fact, the notion of quantum graded A.S.L. is particularly well adapted to
studying both ${\mathcal O}_q(G_{m,n}(\k))$ and its associated quantum Schubert
varieties, which are also quantum graded A.S.L. as shown in [LenRi(2)]. Here,
we use this notion to study quantum Schubert varieties from the point of view
of noncommutative algebraic geometry. Namely, we first show that they are
integral domains which are maximal orders in their division ring of fractions.
Recall that the notion of {\em maximal order} generalises, in the
noncommutative setting, that of a normal domain in commutative algebra. 
Next, we study their regularity properties in the sense of Artin and Schelter
(namely, homological properties such as the Cohen-Macaulay and the Gorenstein
properties).

Also, it turns out that quantum Schubert varieties are strongly linked to
another family of interesting quantum algebras: the quantum determinantal
rings. These are quotients of ${\mathcal O}_q(M_{m,n}(\k))$ by the ideal generated by
quantum minors of a given size. Hence, as a consequence of our results, we are
able to show that quantum determinantal rings are maximal orders in their
division ring of fractions. This generalises results obtained in [LenRi(1)]
where it was shown that quantum determinantal rings are maximal orders under
the restrictive hypotheses that $\k$ is the field of complex numbers and the
deformation parameter $q$ is transcendental over $\Q$. \\

The paper is organised as follows.
Section \ref{sec-basic-definitions} is mainly devoted to recalling basic
definitions and crucial results concerning the quantum algebras we intend to
study.
Section \ref{q-gr-ASL} starts with a short reminder about the notion of
quantum
graded A.S.L. Here, we establish a general criterion
that allows us to prove that, in certain circumstances, a quantum graded
algebra
with a straightening law which is a domain and whose underlying poset has a
single
minimal element is a maximal order in its division ring of fractions. In
Section \ref{q-Schubert-in-G}
we use the results of Section \ref{q-gr-ASL} to show that quantum Schubert
varieties are integral domains and are maximal orders in their
division ring of fractions. Also, we investigate their regularity properties.
It was shown in [LenRi(2)] that quantum Schubert varieties are
AS-Cohen-Macaulay. Here we determine
which of them are AS-Gorenstein. Section \ref{maxor-qdetrings} is devoted to
proving that quantum determinantal rings are maximal orders in their division
ring of fractions. This is shown by using the material of the two preceding
sections.\\

Recall from [McCR; Chap.5, \S 1] that a commutative noetherian domain $A$ is a
maximal order in its quotient field if and only if it is integrally closed.
For this reason, if $A$ is a (noncommutative) domain, then we will say that
$A$ is {\em normal} if it is a maximal order in its division ring of
fractions. This convention differs slightly from classical uses in
noncommutative algebra.

Throughout $\k$ denotes a field. The cardinality of a finite set $X$ is
denoted by $|X|$.

\section{Basic definitions.}\label{sec-basic-definitions}

In this section, we collect some basic definitions and properties about the
objects we intend to study. Most proofs will be omitted since these results
already appear in [KLR] and [LenRi(2)]. Appropriate references will be given
in the text.\\

Let $m,n$ be positive integers. \\

The quantization of the coordinate ring of the affine variety $M_{m,n}(\k)$ of
$m \times n$ matrices with entries in $\k$ is denoted ${\mathcal O}_q(M_{m,n}(\k))$. It
is the $\k$-algebra generated by $mn$ indeterminates $X_{ij}$, with 
$1 \le i \le
m$ and $1 \le j \le n$, subject to the relations:
\[
\begin{array}{ll}  
X_{ij}X_{il}=qX_{il}X_{ij},&\mbox{ for }1\le i \le m,\mbox{ and }1\le j<l\le n\: ;\\ 
X_{ij}X_{kj}=qX_{kj}X_{ij}, & \mbox{ for }1\le i<k \le m, \mbox{ and } 1\le j
\le n \: ; \\  X_{ij}X_{kl}=X_{kl}X_{ij}, & \mbox{ for } 1\le k<i \le m, \mbox{ and }
1\le j<l \le n \: ; \\ X_{ij}X_{kl}-X_{kl}X_{ij}=(q-q^{-1})X_{il}X_{kj}, & 
\mbox{ for } 1\le i<k \le m, \mbox{ and } 1\le j<l \le n. 
\end{array}
\]
To simplify, we write $M_{n}(\k)$ for $M_{n,n}(\k)$. The $m \times n$ matrix
${\bf X}=(X_{ij})$ is called the generic matrix associated with
${\mathcal O}_q(M_{m,n}(\k))$.
\\

As is well known, there exists a $\k$-algebra {\em transpose isomorphism}
between ${\mathcal O}_q(M_{m,n}(\k))$ and ${\mathcal O}_q(M_{n,m}(\k))$, see [LenRi(2); Remark 
3.1.3]. Hence, from now on, we assume that $m \le n$, without loss of
generality.\\

An index pair (in $\{1,\dots,m\}\times\{1,\dots,n\}$) is a pair $(I,J)$ such
that $I \subseteq \{1,\dots,m\}$ and $J \subseteq \{1,\dots,n\}$ are subsets
with the same cardinality. Hence, an index pair is given by an integer $t$
such that $1 \le t \le m$ and ordered sets 
$I=\{i_1 < \dots < i_t\} \subseteq \{1,\dots,m\}$
and $J=\{j_1 < \dots < j_t\} \subseteq \{1,\dots,n\}$. To any such index pair
we associate the quantum minor 
\[ 
[I|J] = \sum_{\sigma\in {\mathfrak S}_t}
(-q)^{\ell(\sigma)} X_{i_{\sigma(1)}j_1} \dots X_{i_{\sigma(t)}j_t} . 
\] 
The
set of all index pairs is denoted by $\Delta_{m,n}$. Since $\Delta_{m,n}$ is
in one-to-one correspondence with the set of all quantum minors of
${\mathcal O}_q(M_{m,n}(\k))$, we will often identify these two sets. The set
$\Delta_{m,n}$ is equipped with the partial order $\le_\st$ defined in
[LenRi(2); Section 3.5]. 
Namely, if $(I,J)$ and $(K,L)$ are index
pairs with $I=\{i_1< \dots <i_u\}, J=\{j_1< \dots <j_u\},  
K=\{k_1< \dots <k_v\}$ and $L=\{l_1< \dots <l_v\}$ 
then 
\[ (I,J) \le_\st (K,L)
\Longleftrightarrow \left\{
\begin{array}{l}
u \ge v, \cr 
i_s \le k_s \quad\mbox{for}\quad 1 \le s \le v , \cr
j_s \le l_s \quad\mbox{for}\quad 1 \le s \le v .
\end{array}
\right.
\]~\\

We now consider the {\it quantization of the coordinate ring of the
grassmannian of $m$-dimensional subspaces of $\k^n$}, denoted by
${\mathcal O}_q(G_{m,n}(\k))$. This is defined to be the subalgebra of
${\mathcal O}_q(M_{m,n}(\k))$ generated by the $m \times m$ quantum minors. \\

An index set (in $\{1,\dots,n\}$) is a subset $I=\{i_1 < \dots < i_m\}
\subseteq \{1,\dots,n\}$. To any index set we associate the maximal quantum
minor $[\{1,\dots,m\}|I]$ of ${\mathcal O}_q(M_{m,n}(\k))$ which is, thus, an element of
${\mathcal O}_q(G_{m,n}(\k))$. The set of all index sets is denoted by $\Pi_{m,n}$.
Since $\Pi_{m,n}$ is in one-to-one correspondence with the set of all maximal
quantum minors of ${\mathcal O}_q(M_{m,n}(\k))$, we will often identify these two sets.
The map $\Pi_{m,n} \longrightarrow \Delta_{m,n}$ given by $I \mapsto
(\{1,\dots,m\},I)$, identifies $\Pi_{m,n}$ with a subset of $\Delta_{m,n}$.
Hence, the partial order $\le_\st$ induces a partial order on $\Pi_{m,n}$ that
we still denote by $\le_\st$. Clearly, if $I=\{i_1< \dots <i_m\}$ and
$J=\{j_1< \dots <j_m\}$ are two index sets, we have \[ I\le_\st J
\Longleftrightarrow i_s \le j_s \quad\mbox{for}\quad 1 \le s \le m. \]

The order in which the $\k$-algebras ${\mathcal O}_q(M_{m,n}(\k))$ and
${\mathcal O}_q(G_{m,n}(\k))$ have been introduced above is forced upon us by the
definition of the quantum grassmannian. Despite this, from our point of view,
the more fundamental object is ${\mathcal O}_q(G_{m,n}(\k))$, rather than
${\mathcal O}_q(M_{m,n}(\k))$ and we concentrate on ${\mathcal O}_q(G_{m,n}(\k))$ in this paper.
Indeed, many desireable properties shared by these two algebras are more
easily proven for ${\mathcal O}_q(G_{m,n}(\k))$. Then, the corresponding property can be
transfered to ${\mathcal O}_q(M_{m,n}(\k))$ using the {\it dehomogenisation map}
$D_{m,n}$, introduced in [KLR], which relates the two algebras. We briefly
recall the definition of this map. Recall from [LenRi(2); \S 3.5] that to any
$(I,J) \in \Delta_{m,n}$ we associate $K_{(I,J)}\in\Pi_{m,m+n}$ where, if
$I=\{i_1<\dots<i_t\}$ and $J=\{j_1<\dots<j_t\}$, for $1 \le t \le m$, we set
$K_{(I,J)}:=\{j_1,\dots,j_t,n+1\dots,n+m\}\setminus\{n+m+1-i_1,\dots,n+m+1-i_t\}$.
Then, letting $M=\{n+1,\dots,n+m\}$, the map $\delta_{m,n} \, : \,
\Delta_{m,n} \longrightarrow \Pi_{m,m+n} \setminus \{M\}$, $(I,J) \mapsto
K_{(I,J)}$ is an isomorphism of partially ordered sets, see [LenRi(2); 3.5.2].
Now, it can be shown that there exists a $\k$-algebra isomorphism 
\[
\begin{array}{ccrcl}
D_{m,n} & : & {\mathcal O}_q(M_{m,n}(\k))[y,y^{-1};\phi] & \longrightarrow &
{\mathcal O}_q(G_{m,m+n}(\k))[[M]^{-1}] 
\end{array}.
\]
given by 
$[I|J]  \mapsto  [K_{(I,J)}][M]^{-1}$ and $y  \mapsto [M]$.
For details about this map see [LenRi(2); \S 3.5].\\

One crucial property in the study of ${\mathcal O}_q(G_{m,n}(\k))$ and
${\mathcal O}_q(M_{m,n}(\k))$ is the existence of {\it standard monomial bases}. By a
{\it standard monomial} in ${\mathcal O}_q(G_{m,n}(\k))$ we mean either $1$ or a product
of the form $[I_1]\dots[I_\ell]$, where $\ell\in\N^\ast$ and $I_1 \le_\st
\dots \le_\st I_\ell \in \Pi_{m,n}$. By [LenRi(2); 3.2.4], the set of standard
monomials is a $\k$-basis of ${\mathcal O}_q(G_{m,n}(\k))$, called the {\it standard
monomial basis} of ${\mathcal O}_q(G_{m,n}(\k))$. A similar notion of standard monomial
can be introduced in ${\mathcal O}_q(M_{m,n}(\k))$ using $\Delta_{m,n}$. Here again, it
turns out that the set of standard monomials is a $\k$-basis of
${\mathcal O}_q(M_{m,n}(\k))$, called the {\it standard monomial basis} of
${\mathcal O}_q(M_{m,n}(\k))$. The latter fact can easily be deduced from the former
using the map $D_{m,n}$, see the proof of [LenRi(2); 3.5.3] for part of the
argument. \\
 
We now introduce the main object of investigation of the present work, namely, quantum analogues of coordinate rings of 
Schubert varieties in the grassmannian. 

\begin{definition} -- \label{def-q-schubert}
Let $\gamma\in\Pi_{m,n}$ and put $\Pi_{m,n}^\gamma=\{\alpha\in\Pi_{m,n} \tq
\alpha\not\ge_\st\gamma\}$. The quantum Schubert variety associated to
$\gamma$ is 
\[ 
{\mathcal O}_q(G_{m,n}(\k))_\gamma := {\mathcal O}_q(G_{m,n}(\k))/\langle
\Pi_{m,n}^\gamma \rangle . 
\]
\end{definition}

\begin{remark} -- \rm \\
(i) Quantum Schubert varieties, as defined above, appear in [LakRe; p.162].
Notice, however, that our conventions differ slightly from those of [LakRe].
However, it is easy to see that the two different conventions produce
isomorphic algebras.\\ (ii) Definition \ref{def-q-schubert} is inspired by the
classical description of the coordinate rings of Schubert varieties in the
grassmannian. For details about this matter, see [GL; \S 6.3.4].
\end{remark}

It turns out that quantum Schubert varieties also have standard monomial
bases. In fact these bases are inherited from the corresponding bases for
${\mathcal O}_q(G_{m,n}(\k))$. One convenient way to show this is by means of the
notion of quantum graded algebra with a straightening law. Hence, we postpone
the details about this point until Section \ref{q-gr-ASL} where this notion is
discussed.

We end this section by showing a technical result that we will use later on.
It is a quantum analogue of Muir's Law of Extensible Minors. In fact, the
result we prove, Proposition \ref{q-Muir} below, is only a special case of the
Quantum Muir's Law of Extensible Minors. For a general result, the reader is
referred to [KroLe; Theorem 3.4]. Even though the result we prove can be
deduced from [KroLe; Theorem 3.4], we have inserted a proof for the
convenience of the reader, since our proof is relatively short. \\

Recall that $n$ is a positive integer and put $F=\{1,\dots,n\}$. 
As is well known, the quantum determinant $[F|F]$ is a central element 
of ${\mathcal O}_q(M_n(\k))$, see, for example, [PW; 4.6.1]. 
Hence, we may form the localisation 
${\mathcal O}_q(GL_n(\k)):=
{\mathcal O}_q(M_n(\k))[[F|F]^{-1}]$. 
By [PW; 5.3.2], ${\mathcal O}_q(GL_n(\k))$ is a Hopf algebra whose 
antipode is the 
anti-automorphism $S$ induced by 
\[
S(X_{ij})  = (-q)^{i-j}[F\setminus\{j\}|F\setminus\{i\}][F|F]^{-1}.
\]
In addition, if $I=\{i_1,\dots,i_t\}$ and $J=\{j_1,\dots,j_t\}$ are 
subsets of $F$ then
\[ 
S([I|J])=(-q)^{(i_1+\dots+i_t)-(j_1+\dots+j_t)}[F \setminus J|F \setminus
I][F|F]^{-1}, 
\] 
see, for example, [KLR; Lemma 4.1].

Now, for $1 \le \ell \le n$, let $\varepsilon_\ell$ denote the element of
$\N^n$ whose only nonzero coordinate equals one and is in the $\ell$-th
position. There is a natural $\N^n \times \N^n$-grading on ${\mathcal
O}_q(M_n)$ relative to which, for $1 \le i,j \le n$, the degree of $X_{ij}$ is
$(\varepsilon_i,\varepsilon_j)$. Clearly, if $I,J$ are subsets of $F$ of the
same cardinality, the degree of $[I|J]$ is $(\sum_{i\in I}
\varepsilon_i,\sum_{j\in J} \varepsilon_j)$.

\begin{proposition} \label{q-Muir} --
Let $P,Q$ be two subsets of $F$ of the same cardinality and denote by
$\overline{P},\overline{Q}$ their respective complements in $F$. Consider
$d\in\N^\ast$ and, for $1 \le s \le d$, elements $c_s \in \k$ and subsets
$I_s,K_s \subseteq P$ and $J_s,L_s \subseteq Q$ such that $|I_s|=|J_s|$ and
$|K_s]=|L_s|$. If the relation $\sum_{s=1}^d c_s[I_s|J_s][K_s|L_s]=0$ holds in
${\mathcal O}_q(M_n(\k))$, then the relation 
\[ 
\sum_{s=1}^d c_s[I_s \cup
\overline{P}|J_s \cup \overline{Q}][K_s \cup \overline{P}|L_s \cup
\overline{Q}]=0 
\] 
holds in ${\mathcal O}_q(M_n(\k))$.
\end{proposition}

\proof We may suppose, without loss of generality, that the products $[I_s|J_s][K_s|L_s]$ in the relation $\sum_{s=1}^d c_s[I_s|J_s][K_s|L_s]=0$
have the same $\N^n \times \N^n$-degree.
Now, let $p$ be the common cardinality of $P$ and $Q$. 
The subalgebra of ${\mathcal O}_q(M_n(\k))$ generated by those $X_{ij}$ such that $i\in P$ and $j\in Q$ is isomorphic to ${\mathcal O}_q(M_p(\k))$. 
Hence, we may consider the relation $\sum_{s=1}^d c_s[I_s|J_s][K_s|L_s]=0$ as a relation in ${\mathcal O}_q(M_p(\k))$ and apply to this relation
the antipode of ${\mathcal O}_q(GL_p(\k))$. This yields the relation 
\[
\sum_{s=1}^d c_s[Q \setminus L_s|P \setminus K_s][Q \setminus J_s|P \setminus I_s]=0
\]
in ${\mathcal O}_q(M_n(\k))$. (Notice that the $\N^n \times \N^n$-homogeniety of the relation has been used here to cancel out the various powers
of $q$ occuring from the application of the antipode.)
Now, applying the antipode of ${\mathcal O}_q(GL_n(\k))$ to this relation gives us the relation
\[
\sum_{s=1}^d c_s[I_s \cup \overline{P}|J_s \cup \overline{Q}][K_s \cup \overline{P}|L_s \cup \overline{Q}]=0
\]
in ${\mathcal O}_q(M_n(\k))$.\qed

\section{Quantum graded algebras with a straightening law.} \label{q-gr-ASL}

In this section, we start reviewing the notion of quantum graded algebra with
a straightening law, as introduced and studied in [LenRi(2)]. Next, we give a
criterion that allows us to show that, under certain hypotheses, a quantum
graded algebra with a straightening law which is a domain and whose underlying
partially 
ordered set has a single minimal element is a maximal order in its division
ring of fractions.

\subsection{Short reminder.} \label{q-gr-ASL-reminder}

In this subsection, we recall the notion of a quantum graded algebra with a
straightening law (on a partially ordered set $\Pi$). 
We also recall various properties
of such algebras that we will use later.\\

Let $A$ be an algebra and $\Pi$ a finite subset of 
elements of $A$ with a partial order $<_\st$. A  
{\em standard monomial} on $\Pi$ is an element
of $A$ which is either $1$ or of the form $\alpha_1\dots\alpha_s$, 
for some $s\geq 1$, where $\alpha_1,\dots,\alpha_s \in \Pi$ and
$\alpha_1\le_\st\dots\le_\st\alpha_s$. 

\begin{subdefinition} \label{recall-q-gr-asl} -- 
Let $A$ be an ${\mathbb N}$-graded $\k$-algebra and $\Pi$ a finite subset 
of $A$ equipped with a partial order $<_\st$. 
We say that $A$ is a {\em quantum graded algebra with a straightening law} 
({\em quantum graded A.S.L.} for short) on the poset $(\Pi,<_\st)$ 
if the following conditions are satisfied.\\
(1) The elements of $\Pi$ are homogeneous with positive degree.\\
(2) The elements of $\Pi$ generate $A$ as a $\k$-algebra.\\
(3) The set of standard monomials on $\Pi$ is a linearly independent set.\\
(4) If $\alpha,\beta\in\Pi$ are not comparable for $<_\st$, 
then $\alpha\beta$ 
is a linear combination of terms $\lambda$ or $\lambda\mu$, where 
$\lambda,\mu\in\Pi$, $\lambda\le_\st\mu$ and $\lambda<_\st\alpha,\beta$.\\
(5) For all $\alpha,\beta\in\Pi$, there exists $c_{\alpha\beta} \in \k^\ast$ 
such that $\alpha\beta-c_{\alpha\beta}\beta\alpha$ is a linear combination of 
terms $\lambda$ or $\lambda\mu$, where $\lambda,\mu\in\Pi$,
$\lambda\le_\st\mu$ and $\lambda<_\st\alpha,\beta$.
\end{subdefinition}

By [LenRi(2); Proposition 1.1.4], if $A$ is a quantum graded A.S.L. on the
partially 
ordered set $(\Pi,<_\st)$, then the set of standard monomials on $\Pi$ forms a
$\k$-basis of $A$. Hence, in the presence of a standard monomial basis, the
structure of a quantum graded A.S.L. may be seen as providing substancial
further information on the way standard monomials multiply and commute.

\begin{subexample} -- \rm
As is well known, the algebra ${\mathcal O}_q(M_{m,n}(\k))$ is $\N$-graded, by putting
the canonical generators in degree one. Now, since ${\mathcal O}_q(G_{m,n}(\k))$ is a
subalgebra of ${\mathcal O}_q(M_{m,n}(\k))$ generated by homogeneous elements,
${\mathcal O}_q(G_{m,n}(\k))$ inherits a natural $\N$-grading from that of
${\mathcal O}_q(M_{m,n}(\k))$. In fact, beyond the existence of standard monomial bases
for ${\mathcal O}_q(G_{m,n}(\k))$ and ${\mathcal O}_q(M_{m,n}(\k))$, as mentioned in Section
\ref{sec-basic-definitions}, we have that ${\mathcal O}_q(G_{m,n}(\k))$ is a quantum
graded algebra on $(\Pi_{m,n},\le_\st)$ and that ${\mathcal O}_q(M_{m,n}(\k))$ is a
quantum graded algebra on $(\Delta_{m,n},\le_\st)$, see [LenRi(2); Theorem
3.4.4 and 3.5.3].
\end{subexample}

From our point of view, one important feature of quantum graded A.S.L. 
is the
following. Let $A$ be a $\k$-algebra which is a quantum graded A.S.L. on the
set $(\Pi,\le_\st)$. A subset $\Omega$ of $\Pi$ will be called a $\Pi$-ideal
if it is an ideal of the partially 
ordered set $(\Pi,\le_\st)$ in the sense of lattice
theory; that is, if it satisfies the following property: if $\alpha\in\Omega$
and if $\beta\in\Pi$, with $\beta \le_\st \alpha$, then $\beta\in\Omega$. We
can consider the quotient $A/\langle\Omega\rangle$ of $A$ by the ideal
generated by $\Omega$. Clearly, it is still a graded algebra and it is
generated by the images in $A/\langle\Omega\rangle$ of the elements of
$\Pi\setminus\Omega$. The important point here is that
$A/\langle\Omega\rangle$ inherits from $A$ 
a natural quantum graded A.S.L. structure on
$\Pi\setminus\Omega$ (or, more precisely, 
on the canonical image of $\Pi\setminus\Omega$
in $A/\langle\Omega\rangle$). In particular, the set of
homomorphic images in $A/\langle\Omega\rangle$ of the standard monomials of
$A$ which either equal $1$ or are of the form $\alpha_1 \dots \alpha_t$
($t\in\N^\ast$) and $\alpha_1\notin\Omega$ form a $\k$-basis for
$A/\langle\Omega\rangle$. The reader will find all the necessary details in \S
1.2 of [LenRi(2)].

\begin{subexample} -- \rm \label{example-q-gr-ASL}
Let $\gamma \in\Pi_{m,n}$. It is clear that the set $\Pi_{m,n}^\gamma$
introduced in Definition \ref{def-q-schubert} is a $\Pi_{m,n}$-ideal. Hence,
the discussion above shows that the quantum Schubert variety
${\mathcal O}_q(G_{m,n}(\k))_\gamma$ is a quantum graded A.S.L. on the canonical image
in ${\mathcal O}_q(G_{m,n}(\k))_\gamma$ of $\Pi_{m,n} \setminus \Pi_{m,n}^\gamma$. In
particular, the canonical images in ${\mathcal O}_q(G_{m,n}(\k))_\gamma$ of the standard
monomials of ${\mathcal O}_q(G_{m,n}(\k))$ which either equal to $1$ or are of the form
$[I_1]\dots[I_t]$, for some $t\geq 1$ and with $\gamma \le_\st [I_1]$, form a
$\k$-basis of ${\mathcal O}_q(G_{m,n}(\k))_\gamma$.\\
\end{subexample}

\begin{subremark} -- \rm \label{remark-single-min-elt}
Let $\gamma \in\Pi_{m,n}$. As mentioned in Example \ref{example-q-gr-ASL}, the
quantum Schubert variety ${\mathcal O}_q(G_{m,n}(\k))_\gamma$ is a quantum graded A.S.L.
on the canonical image in ${\mathcal O}_q(G_{m,n}(\k))_\gamma$ of $\Pi_{m,n} \setminus
\Pi_{m,n}^\gamma$. At this point, it is worth noting that the set $\Pi_{m,n}
\setminus \Pi_{m,n}^\gamma$ has a single minimal element, namely $\gamma$.
\end{subremark}

We end this subsection 
by recalling from [LenRi(2); Proposition 1.1.5] 
the formula which gives the Gelfand-Kirillov dimension of a 
quantum graded A.S.L. Recall
that, if $(\Pi,\le_\st)$ is a partially 
ordered set, the rank of an element $\pi\in\Pi$, denoted $\rk\pi$, is the greatest integer 
$k\in\N$ such that there exists a chain 
$\pi_1 <_\st \dots <_\st \pi_{k-1} <_\st \pi_k=\pi$ of elements of $\Pi$. Then, we define the rank of $\Pi$ by $\rk\Pi=\max\{\rk\pi,\,\pi\in\Pi\}$. 
Then, we have the following proposition.

\begin{subproposition} -- \label{GKdim-q-gr-ASL} 
Let $A$ be a quantum graded A.S.L. on $(\Pi,\le_\st)$; then $\GKdim A =\rk\Pi$.
\end{subproposition}

\begin{subcorollary} -- \label{GKdim-q-schubert-grass}
Let $\gamma = (\gamma_1, \dots, \gamma_m)\in\Pi_{m,n}$. 
Then
\[
\GKdim {\mathcal O}_q(G_{m,n}(\k))_\gamma=m(n-m)+\frac{m(m+1)}{2}-
\left(\sum_{i=i}^m\gamma_i\right)
+1.
\]
\end{subcorollary}

\proof It is well known that $\rk (\Pi_{m,n}\setminus\Pi_{m,n}^\gamma) =
m(n-m)+\frac{m(m+1)}{2}-
\left(\sum_{i=i}^m\gamma_i\right) +1$, see [BV; 5.12]; so, 
the result follows from  Proposition \ref{GKdim-q-gr-ASL}.\qed

\subsection{Quantum graded A.S.L. and the maximal order
property.}\label{par-normality}
In this subsection we are interested in quantum graded A.S.L. whose associated
poset has a single minimal element. Hence, let $A$ be a quantum graded A.S.L.
on the poset $(\Pi,\le_\st)$, and assume that $(\Pi,\le_\st)$ has a single
minimal element, denoted $\gamma$.
We know 
that $\gamma$ is a regular normal element
of $A$, by [LenRi(2); Lemma 1.2.1]. 
Hence, we may form the localisation, $A[\gamma^{-1}]$, of $A$
with respect to the powers of $\gamma$ and the canonical map $A
\longrightarrow A[\gamma^{-1}]$ is injective. Notice that
quantum Schubert varieties are examples of such
algebras, as mentioned in
Remark
\ref{remark-single-min-elt}.\\

Our first interest is in studying the ideal $\langle\gamma\rangle$ of $A$. \\

To each element $\sigma\in\Pi$, we associate the subset $\Pi^\sigma$ 
of $\Pi$ defined by
\[
\Pi^\sigma = \{\pi\in\Pi \tq \pi\not\geq_\st\sigma\}.
\]
It is clear that $\Pi^\sigma$ is a $\Pi$-ideal. 
In addition, we let $I_\sigma$ be the ideal of $A$ generated by $\Pi^\sigma$:
\[
I_\sigma=\langle\Pi^\sigma\rangle.
\]
It is clear that, if $\sigma$ and $\tau$ are elements of $\Pi$ such that $\sigma\le_\st\tau$, then $\Pi^\sigma \subseteq \Pi^\tau$ and hence 
$I_\sigma \subseteq I_\tau$. Finally, a last piece of notation: let $\sigma\in\Pi$; an element $\tau$ is called an upper neighbour of $\sigma$ 
if $\sigma <_\st\tau$ and there is no element $\nu\in\Pi$ such that  $\sigma <_\st \nu <_\st \tau$. Clearly, for any $\sigma\in\Pi$, the set of 
upper neighbours of $\sigma$ is a (finite) subset of $\Pi$ which is empty if and only if $\sigma$ is maximal.

\begin{sublemma} \label{abstract-pieri} -- 
We keep the notation introduced above.\\
(i) Let $\psi$ be the automorphism of $A$ associated to the regular normal 
element $\gamma$; that is, $\gamma a=\psi(a)\gamma$, for all $a\in A$. Then
$\psi(I_\tau)=I_\tau$ for all $\tau\in\Pi$.\\ 
(ii) If $\{\gamma\}\subsetneqq \Pi$, then 
the following formula holds, where the intersection is taken over all the
upper neigbours of $\gamma$:

\[
\langle\gamma\rangle=\bigcap I_\tau.
\]
\end{sublemma}

\proof (i) By condition (5) of Definition \ref{recall-q-gr-asl}, the element 
$\gamma$
commutes up to a non zero scalar with each element of $\Pi$. Thus, 
each $\tau\in\Pi$ is an eigenvector (with non zero eigenvalue) of $\psi$; so,
the statement is clear.\\ 
(ii) As $\gamma \in \Pi^\tau$, for all upper
neighbour $\tau$ of $\gamma$, the inclusion
$\langle\gamma\rangle\subseteq\bigcap I_\tau$ is clear. 
Let us now obtain the reverse
inclusion. First, notice that any element of $\Pi$ different from $\gamma$
must be greater than or equal to some upper neighbour of $\gamma$; this is an
easy consequence of the fact that $\Pi$ is a finite partially 
ordered set whose unique
minimal element is $\gamma$. On the other hand, for all $\tau\in\Pi$, since
$\Pi^\tau$ is a $\Pi$-ideal, the ideal $I_\tau$ is the vector space generated
by standard monomials involving an element of $\Pi^\tau$, see [LenRi(2);
Proposition 1.2.5]. Hence, $I_\tau$ is the vector space generated by standard
monomials of the form $\alpha_1 \dots \alpha_r$, with 
$\alpha_1 \le_\st
\dots \le_\st \alpha_r \in\Pi$ and such that $\alpha_1\not\geq_\st\tau$. Since, in
addition, the standard monomials form a basis of $A$, it follows that $\bigcap
I_\tau$ (where the intersection is taken over all the upper neighbours of
$\gamma$) is the vector space generated by standard monomials of the form
$\alpha_1 \dots \alpha_r$, with $\alpha_1 \le_\st \dots \le_\st \alpha_r
\in\Pi$ and such that $\alpha_1$ is not greater than or equal to any upper
neighbour of $\gamma$. By the above comment, this forces $\alpha_1
= \gamma$. 
The inclusion $\bigcap I_\tau\subseteq\langle\gamma\rangle$ now follows.\qed

\begin{subproposition} -- \label{CS-normality-q-gr-ASL}
We keep the notation introduced above. Assume that $A$ is a domain such that 
$I_\tau$ is a completely prime
ideal of $A$ for any upper neighbour $\tau$ of $\gamma$,
and that $A[\gamma^{-1}]$ is a maximal order in its division ring
of fractions. Then $A$ is a maximal order in its division ring of fractions.
\end{subproposition}

\proof First, recall from [LenRi(2); Lemma 1.2.3] that $A$ is noetherian.
Notice, in addition, that if $\Pi=\{\gamma\}$, then $A$ is a commutative
polynomial ring in one indeterminate; so that $A$ is clearly a maximal order
in its division ring of fractions. 
 Now, assume that $\{\gamma\} \subsetneqq \Pi$.  
By Lemma
\ref{abstract-pieri} and the hypotheses made on $A$, we are in position to
apply [R; Lemma 1.1] which gives the result.\qed

\section{Quantum Schubert varieties in the grassmannian.}
\label{q-Schubert-in-G}
Let $m,n$ be positive integers. \\

As discussed above, to each $\gamma\in\Pi_{m,n}$, we may associate the
$\k$-algebra ${\mathcal O}_q(G_{m,n}(\k))_\gamma$, which is a quantum deformation of the
coordinate ring of a Schubert variety. The aim of this section is to study
these rings. In the first subsection, we will show that they are normal
integral domains. In the second section, we will study them from the point of
view of regularity conditions.

\subsection{Integrality and normality.}\label{int-and-norm}

Our aim in this subsection is to prove that the quantum Schubert variety
${\mathcal O}_q(G_{m,n}(\k))_\gamma$ is a normal domain for any $\gamma\in\Pi_{m,n}$.
Here, by {\em normal domain}, we mean an integral domain which is a maximal
order in its division ring of fractions. To achieve this goal, we will be
naturally led to use certain subalgebras of ${\mathcal O}_q(M_{m,n}(\k))$ that we now
define. \\

Let us start by introducing some convenient notation. To each
$\gamma=(\gamma_1,\dots,\gamma_m)\in\Pi_{m,n}$, with $1 \le \gamma_1< \dots <
\gamma_m\le n$, we associate the substet ${\mathcal L}_\gamma$ of $\{1,\dots,m\} \times
\{1,\dots,n\}$ defined by 
\[ 
{\mathcal L}_\gamma = \{(i,j) \in \{1,\dots,m\} \times
\{1,\dots,n\} \;\tq\; j > \gamma_{m+1-i} \quad\mbox{and}\quad j\neq \gamma_\ell
\quad\mbox{for}\quad 1 \le \ell \le m\}, 
\] 
which we call the {\em ladder}
associated with $\gamma$.

It follows from the definition of ${\mathcal L}_\gamma$ that for each
$(i,j)\in{\mathcal L}_\gamma$ the set 
$\{\gamma_1,\dots,\gamma_m\} \setminus\{\gamma_{m+1-i}\} \cup \{j\}$ is a
subset of $\{1,\dots,n\}$ containing $m$ distinct elements. 
Hence, it makes sense to
associate to this subset the maximal quantum minor 
$m_{ij}:=[\{\gamma_1,\dots,\gamma_m\} \setminus\{\gamma_{m+1-i}\} 
\cup \{j\}]$ 
of ${\mathcal O}_q(M_{m,n}(\k))$. We then set 

\[
\M_\gamma=\{m_{ij} \in \Pi_{m,n} \tq (i,j)\in{\mathcal L}_\gamma\} \subseteq \Pi_{m,n}.
\]

\begin{subremark} -- \label{rem-on-mij} \rm 
Let $\gamma=(\gamma_1,\dots,\gamma_m)\in\Pi_{m,n}$.\\ (i) Consider
$j,l\in\{1,\dots,m\}$ and suppose $j \notin\{\gamma_1,\dots,\gamma_m\}$. Then,
clearly, $j < \gamma_l$ implies that
$[\{\gamma_1,\dots,\gamma_m\}\setminus\{\gamma_l\}\cup\{j\}] <_\st \gamma$
while $j > \gamma_l$ implies that
$[\{\gamma_1,\dots,\gamma_m\}\setminus\{\gamma_l\}\cup\{j\}] >_\st \gamma$.\\
(ii) Hence, the elements of $\M_\gamma$ are nothing but the elements of
$\Pi_{m,n}$ that are greater than $\gamma$ with respect to the partial 
order $\le_\st$
and differ from $\gamma$ by exactly one column index.
\end{subremark}

\begin{subdefinition} -- \label{def-q-ladder-matrix} 
Let $\gamma=(\gamma_1,\dots,\gamma_m)\in\Pi_{m,n}$, with $1 \le \gamma_1< \dots < \gamma_m\le n$. The quantum ladder 
matrix ring associated with $\gamma$, denoted ${\mathcal O}_q(M_{m,n,\gamma}(\k))$, is the $\k$-subalgebra of ${\mathcal O}_q(M_{m,n}(\k))$ generated by the elements 
$X_{ij}\in {\mathcal O}_q(M_{m,n}(\k))$ such that $(i,j)\in{\mathcal L}_\gamma$.
\end{subdefinition}

Let us discuss an example to clarify the definition.

\begin{subexample} -- \rm
We put $(m,n)=(3,7)$ and
$\gamma=(\gamma_1,\gamma_2,\gamma_3)=(1,3,6)\in\Pi_{3,7}$. In the $3 \times 7$
generic matrix ${\bf X}=\left(X_{ij}\right)$ associated to
${\mathcal O}_q(M_{3,7}(\k))$, put a bullet on each row as follows: on the first row,
the bullet is in column $6$ because $\gamma_3$ is $6$, on the second row, the
bullet is in column $3$ because $\gamma_2$ is $3$ and on the third row, 
the bullet is
in column $1$ because $\gamma_1 = 1$. 
Now, in each position which is to the left of a bullet, or
which is below a bullet, put a star. To finish, place 
$X_{ij}$ in any  position $(i,j)$ that has not yet been filled. 
We obtain 
\[
\left(
\begin{array}{ccccccc}
 \ast & \ast & \ast & \ast & \ast & \bullet & X_{17} \cr   
 & & & & & & \cr 
\ast & \ast & \bullet & X_{24} & X_{25} & \ast & X_{27} \cr
 & & & & & & \cr
\bullet & X_{32} & \ast & X_{34} & X_{35} & \ast & X_{37} \cr
\end{array}
\right) .
\]
By definition, the ladder quantum matrix ring associated to $\gamma=(1,3,6)$ is the subalgebra of ${\mathcal O}_q(M_{3,7}(\k))$ generated by the elements
$X_{17}, X_{24}, X_{25}, X_{27}, X_{32}, X_{34}, X_{35}, X_{37}$.  
\end{subexample}

Our aim now is to show that 
the localisation of
${\mathcal O}_q(G_{m,n}(\k))_\gamma$ at the powers of the image of $\gamma$ in
${\mathcal O}_q(G_{m,n}(\k))_\gamma$ 
is isomorphic to a skew Laurent extension of
${\mathcal O}_q(M_{m,n,\gamma}(\k))$.

\begin{sublemma} -- \label{relations-ladder-minors}
Let $\gamma=(\gamma_1,\dots,\gamma_m)\in\Pi_{m,n}$, with $1 \le \gamma_1< \dots < \gamma_m\le n$. 
For $(i,j),(k,l)\in{\mathcal L}_\gamma$, 
the following relations hold in ${\mathcal O}_q(G_{m,n}(\k))$:\\
(i) if $i=k$ and $j<l$, then $m_{ij}m_{kl}=qm_{kl}m_{ij}$;\\
(ii) if $i<k$ and $j=l$, then $m_{ij}m_{kl}=qm_{kl}m_{ij}$; \\
(iii) if $i<k$ and $j>l$, then $m_{ij}m_{kl}=m_{kl}m_{ij}$;\\
(iv) if $i<k$ and $j<l$, then $m_{ij}m_{kl}-m_{kl}m_{ij}=(q-q^{-1})m_{il}m_{kj}$;\\
(v) $\gamma m_{ij}=qm_{ij}\gamma$.
\end{sublemma}

\proof The proof is an easy application of Proposition \ref{q-Muir}. We give
details for (iv) leaving the other (easier) cases to the reader.

First, note that, if we set $R=\{\gamma_1,\dots,\gamma_m\}
\setminus\{\gamma_{m+1-i},\gamma_{m+1-k}\}$, then 
$m_{ij}=[\{\gamma_1,\dots,\gamma_m\} \setminus\{\gamma_{m+1-i}\} \cup
\{j\}]=[R \cup \{\gamma_{m+1-k},j\}]$ and $m_{kl}=[\{\gamma_1,\dots,\gamma_m\}
\setminus\{\gamma_{m+1-k}\} \cup \{l\}]=[R \cup \{\gamma_{m+1-i},l\}]$. 
In addition, $\gamma_{m+1-i}<j$, 
since $(i,j)$ and $(k,l)$ are in ${\mathcal L}_\gamma$. 
Hence, 
\[
\gamma_{m+1-k} < \gamma_{m+1-i} < j < l.
\]
It follows that, in ${\mathcal O}_q(M_n(\k))$, we have the relation 
\[
\begin{array}{c}
[1,2|\gamma_{m+1-k},j][1,2|\gamma_{m+1-i},l]-[1,2|\gamma_{m+1-i},l][1,2|\gamma_{m+1-k},j]\cr
\cr
=(q-q^{-1})[1,2|\gamma_{m+1-k},l][1,2|\gamma_{m+1-i},j].
\end{array}
\] 
(This is an immediate consequence of the relation $[13][24]-[24][13]=(q-q^{-1})[14][23]$ that holds in ${\mathcal O}_q(G_{2,4}(\k))$, see the introduction of [KLR],
using a suitable injection
of ${\mathcal O}_q(G_{2,4}(\k))$ in ${\mathcal O}_q(M_n(\k))$.) Now, applying Proposition \ref{q-Muir} to this relation, with $\overline{P}=\{3,\dots,m\}$ and 
$\overline{Q}=R$ gives the relation 
\[
\begin{array}{c}
[1,\dots,m|R \cup \{\gamma_{m+1-k},j\}][1,\dots,m|R \cup\{\gamma_{m+1-i},l\}] \cr\cr
\qquad -[1,\dots,m|R \cup\{\gamma_{m+1-i},l\}][1,\dots,m|\{R \cup\gamma_{m+1-k},j\}]\cr
\cr
=(q-q^{-1})[1,\dots,m|R \cup\{\gamma_{m+1-k},l\}][1,\dots,m|R \cup\{\gamma_{m+1-i},j\}]. 
\end{array}
\]
If we view ${\mathcal O}_q(G_{m,n}(\k)$ as the subalgebra of ${\mathcal O}_q(M_n(\k))$ 
generated by the $m \times m$ minors built on the first $m$ rows of the generic
matrix of ${\mathcal O}_q(M_n(\k))$, this gives the required relation. \qed

\begin{subremark} --\label{remark-on-ladder} \rm \\
(i) In view of the defining relations of ${\mathcal O}_q(M_{m,n}(\k))$, it is clear that
there exists a $\k$-algebra automorphism $\psi \, : \, {\mathcal O}_q(M_{m,n}(\k))
\longrightarrow {\mathcal O}_q(M_{m,n}(\k))$ such that $\psi(X_{ij})=qX_{ij}$, 
for each $i,j$.\\
(ii) Let $\gamma\in\Pi_{m,n}$. It is not difficult to check that the quantum
ladder matrix ring, ${\mathcal O}_q(M_{m,n,\gamma}(\k))$, is isomorphic to an iterated
skew polynomial extension of $\k$ obtained by inserting the generators
$X_{ij}$, with $(i,j)\in{\mathcal L}_\gamma$, in lexicographic order. In fact,
${\mathcal O}_q(M_{m,n,\gamma}(\k))$ is isomorphic to the $\k$-algebra generated by
indeterminates $X_{ij}$, with $(i,j)\in{\mathcal L}_\gamma$, subject to the relations
imposed by the fact that  $X_{ij} \in {\mathcal O}_q(M_{m,n}(\k))$.
In addition $\psi$
clearly restricts to a $\k$-algebra automorphism of
${\mathcal O}_q(M_{m,n,\gamma}(\k))$, that we also denote by $\psi$.\\ 
(iii) Let
$\gamma\in\Pi_{m,n}$. By point (ii) above and standard results on
Gelfand-Kirillov dimension, one has $\GKdim
{\mathcal O}_q(M_{m,n,\gamma}(\k))=|{\mathcal L}_\gamma|$. Now, clearly, $\psi$ is a locally
algebraic automorphism of ${\mathcal O}_q(M_{m,n,\gamma}(\k))$ in the sense of [KL; \S
12.3]. It follows, see [KL; \S 12.3], that:
\[
\GKdim {\mathcal O}_q(M_{m,n,\gamma}(\k))[Y,Y^{-1};\psi]=|{\mathcal L}_\gamma|+1=
m(n-m)+\frac{m(m+1)}{2}-\left(\sum_{i=i}^m\gamma_i\right) +1.
\]
\end{subremark}

Let $\gamma\in\Pi_{m,n}$. The homomorphic image in ${\mathcal O}_q(G_{m,n}(\k))_\gamma$
of an element $x\in {\mathcal O}_q(G_{m,n}(\k))$ will be denoted $\overline{x}$. Theorem
\ref{dehom-q-schubert} establishes a strong link between the localisation of
${\mathcal O}_q(G_{m,n}(\k))_\gamma$ at powers of $\overline{\gamma}$ and a skew Laurent
extension of the quantum ladder matrix ring ${\mathcal O}_q(M_{m,n,\gamma}(\k))$. Recall
that ${\mathcal O}_q(G_{m,n}(\k))_\gamma$ is a quantum graded A.S.L. on the set
$\Pi_{m,n}\setminus\Pi_{m,n}^\gamma$ (identified with its image in
${\mathcal O}_q(G_{m,n}(\k))_\gamma$). 
The element
$\overline{\gamma}$ is the unique minimal element of
$\Pi_{m,n}\setminus\Pi_{m,n}^\gamma$; and so
$\overline{\gamma}$ is a regular normal
element of ${\mathcal O}_q(G_{m,n}(\k))_\gamma$.
Thus, we may form the localisation
${\mathcal O}_q(G_{m,n}(\k))_\gamma[\overline{\gamma}^{-1}]$ of
${\mathcal O}_q(G_{m,n}(\k))_\gamma$ at the powers of $\overline{\gamma}$, as indicated
in the
introduction to Subsection \ref{par-normality}.

\begin{subtheorem} -- \label{dehom-q-schubert}
Let $\gamma\in\Pi_{m,n}$. There exists a $\k$-algebra isomorphism 
\[
\begin{array}{ccrcl}
d_\gamma & : & {\mathcal O}_q(M_{m,n,\gamma}(\k))[Y,Y^{-1};\phi] & \longrightarrow &
{\mathcal O}_q(G_{m,n}(\k))_\gamma\,[\overline{\gamma}^{-1}] 
\end{array} .
\]
sending $
X_{ij}$ to $\overline{m_{ij}}$ and $Y$ to $\overline{\gamma}$. 
\end{subtheorem}

\proof The existence of the $\k$-algebra morphism $d_\gamma$ is clear from
Lemma \ref{relations-ladder-minors} and Remark \ref{remark-on-ladder}. We
now show that $d_\gamma$ is onto, by showing that
${\mathcal O}_q(G_{m,n}(\k))_\gamma[\overline{\gamma}^{-1}]$ is generated as a
$\k$-algebra by $\overline{\gamma}$, $\overline{\gamma}^{-1}$ and
$\overline{m_{ij}}$, $(i,j)\in{\mathcal L}_\gamma$. Clearly,
${\mathcal O}_q(G_{m,n}(\k))_\gamma[\overline{\gamma}^{-1}]$ is generated as a
$\k$-algebra by $\overline{\gamma}$, $\overline{\gamma}^{-1}$ and the images
in ${\mathcal O}_q(G_{m,n}(\k))_\gamma$ of the minors in
$\Pi_{m,n}\setminus\Pi_{m,n}^\gamma$; so, what we must show is that
$\overline{x} \in\im d_\gamma$, 
for all
$x\in\Pi_{m,n}\setminus\Pi_{m,n}^\gamma$. To
each $x = [j_1,\dots,j_m]\in \Pi_{m,n}\setminus\Pi_{m,n}^\gamma$, we associate
the number $n(x)$ of elements $j_s$, with $1 \le s \le m$, such that
$j_s\notin\{\gamma_1,\dots,\gamma_m\}$. As mentioned in Remark
\ref{rem-on-mij}, if $x\in\Pi_{m,n}\setminus\Pi_{m,n}^\gamma$ is such that
$n(x) \le 1$, then 
either $\overline{x} = \overline{\gamma}$ or 
$\overline{x} = \overline{m_{ij}}$ for some $ (i,j)\in{\mathcal L}_\gamma$; so that 
$\overline{x} \in \im d_\gamma$. Now, suppose $t$ is an integer in
$\{1,\dots,m-1\}$ with the property that any
$x\in\Pi_{m,n}\setminus\Pi_{m,n}^\gamma$ such that $n(x) \le t$ satisfies
$\overline{x} \in \im d_\gamma$. Consider
$x=[j_1,\dots,j_m]\in\Pi_{m,n}\setminus\Pi_{m,n}^\gamma$ such that $n(x)=t+1$.
In addition, let $1 \le \ell \le m$ be such that
$j_\ell\notin\{\gamma_1,\dots,\gamma_m\}$. The generalised quantum Pl\"ucker
relations of [KLR; Theorem 2.1], applied with $J_1=\emptyset$,
$J_2=\{j_1,\dots,j_m\}\setminus\{j_\ell\}$ and
$K=\{\gamma_1,\dots,\gamma_m\}\cup\{j_\ell\}$ gives a relation \[
\sum_{K'\sqcup K''=K} (-q)^{\bullet}[K'][K'' \sqcup J_2]=0 \] in
${\mathcal O}_q(G_{m,n}(\k))$. (Here, by a symbol $(-q)^{\bullet}$, we mean some power
of $-q$ with exponent in $\Z$.) Let us now consider the various terms
$[K'][K'' \sqcup J_2]$ of the above equation. When $K''=\{j_\ell\}$, then
$[K'][K'' \sqcup J_2]=\gamma x$. Otherwise, $K''=\{\gamma_k\}$ for some $1 \le
k \le m$ such that $\gamma_k \notin\{j_1,\dots,j_m\}$. In this case, $[K'][K''
\sqcup
J_2]=[\{\gamma_1,\dots,\gamma_m\}\setminus\{\gamma_k\}\cup\{j_\ell\}][\{j_1,\dots,j_m\}\setminus\{j_\ell\}\cup\{\gamma_k\}]$
(notice that the image of such a term in ${\mathcal O}_q(G_{m,n}(\k))_\gamma$ might very
well be zero). Hence, taking the image of the above relation in
${\mathcal O}_q(G_{m,n}(\k))_\gamma$, we get a relation of the form

\[
\overline{\gamma} \overline{x} = \sum (-q)^{\bullet} \overline{y}\overline{z}
\]
where the sum extends over pairs $(y,z)$ of elements in
$\Pi_{m,n}\setminus\Pi_{m,n}^\gamma$ such that $n(y)=1$ and $n(z)=t$.
Hence,
each 
$\overline{x} \in \im(d_\gamma)$, 
by the induction hypothesis since $\overline{\gamma}$ 
is invertible in ${\mathcal O}_q(G_{m,n}(\k))_\gamma\,[\overline{\gamma}^{-1}]$. 
This
shows that $d_\gamma$ is surjective. Recall from Remark \ref{remark-on-ladder}
that
$\GKdim {\mathcal O}_q(M_{m,n,\gamma}(\k))[Y,Y^{-1};\psi]=
m(n-m)+\frac{m(m+1)}{2}-(\sum_{i=1}^m\gamma_i) +1$.
On the other hand, $\overline{\gamma}$ is a local normal element in the sense
of [KL; \S 12.4]. Hence, by [KL; Theorem 12.4.4] and
Corollary \ref{GKdim-q-schubert-grass}, $\GKdim
{\mathcal O}_q(G_{m,n}(\k))_\gamma[\overline{\gamma}^{-1}]
=m(n-m)+\frac{m(m+1)}{2}-(\sum_{i=1}^m\gamma_i)
+1$.
As ${\mathcal O}_q(M_{m,n,\gamma}(\k))[Y,Y^{-1};\psi]$ is an
integral domain and $d_\gamma$ is surjective, 
the
injectivity of $d_\gamma$ follows, by [KL; Proposition 3.15],.
\qed

\begin{subcorollary} -- \label{q-schubert-normal-domain}
Let $\gamma\in\Pi_{m,n}$. The $\k$-algebra ${\mathcal O}_q(G_{m,n}(\k))_\gamma$ is a normal domain.
\end{subcorollary}

\proof Consider $\tau\in\Pi_{m,n}$. Theorem \ref{dehom-q-schubert} asserts
that there exists a $\k$-algebra isomorphism
${\mathcal O}_q(M_{m,n,\tau}(\k))[Y,Y^{-1};\psi] \cong
{\mathcal O}_q(G_{m,n}(\k))_\tau[\overline{\tau}^{-1}]$. Hence
${\mathcal O}_q(G_{m,n}(\k))_\tau[\overline{\tau}^{-1}]$ is an integral domain. As we
already mentioned, $\overline{\tau}$ is a regular element of
${\mathcal O}_q(G_{m,n}(\k))_\tau$; so the canonical map
${\mathcal O}_q(G_{m,n}(\k))_\tau\longrightarrow{\mathcal O}_q(G_{m,n}(\k))_\tau[\overline{\tau}^{-1}]$
is an injection. It follows that ${\mathcal O}_q(G_{m,n}(\k))_\tau$ is an integral
domain. Hence,
we have proved that any quantum Schubert variety is an integral domain.

It follows that, in the notation of Subsection \ref{par-normality}, the ideal
$I_\tau$ of ${\mathcal O}_q(G_{m,n}(\k))_\gamma$ is a completely prime ideal for all
$\tau\in\Pi_{m,n}\setminus \Pi_{m,n}^\gamma$. Moreover, since
${\mathcal O}_q(G_{m,n}(\k))_\gamma[\overline{\gamma}^{-1}]$ is isomorphic to a
localisation of an iterated Ore extension of $\k$, it is a normal domain by
[MR; V. Proposition 2.5, IV Proposition 2.1]. Hence, Proposition
\ref{CS-normality-q-gr-ASL} applies to the quantum graded A.S.L.
${\mathcal O}_q(G_{m,n}(\k))_\gamma$ (whose underlying poset has a single minimal
element, as noticed in Remark \ref{remark-single-min-elt}); so we conclude
that ${\mathcal O}_q(G_{m,n}(\k))_\gamma$ is a normal domain.\qed

\subsection{The AS-Cohen-Macaulay and AS-Gorenstein properties.}

The following result is Theorem 4.2 of [LenRi(2)].

\begin{subtheorem} -- \label{q-sch-grass-CM} 
Let $\gamma \in \Pi_{m,n}$. The quantum Schubert variety ${\mathcal O}_q(G_{m,n}(\k))_\gamma$ is AS-Cohen-Macaulay.
\end{subtheorem}
 
It is now easy to determine, among quantum Schubert varieties those which are
AS-Gorenstein. Let $A$ be a noetherian $\N$-graded connected $\k$-algebra. For
the definition of the AS-Gorenstein condition for $A$ see [LenRi(2);
Subsection 2.1]. Suppose in addition that $A$ has enough normal elements in
the sense of Zhang, see [LenRi(2); Definition 2.1.3]. Then, $A$ is
AS-Gorenstein if and only if it has finite injective dimension on both sides. 
In particular, if $A$ is commutative, then $A$ is AS-Gorenstein if and only if
it is Gorenstein in the usual sense. For details on these statements, see
Subsection 2.1 of [LenRi(2)] and in particular [LenRi(2); Remark 2.1.10].\\

We need to introduce some more notation. Let
$\gamma=(\gamma_1,\dots,\gamma_m)\in\Pi_{m,n}$. Following section B in Chapter
6 of [BV], we denote by $\beta_0,\dots, \beta_s$ the blocks of consecutive
integers in $\gamma$ and by $\chi_0,\dots,\chi_s$ the gaps between these
blocks. Notice that $\chi_s$ is empty if and only if $\gamma_m=n$.

\begin{subtheorem} -- \label{q-sch-grass-ASG}
Let $\gamma=(\gamma_1,\dots,\gamma_m)\in\Pi_{m,n}$. In the previous notation,
we put $t=s$ if $\gamma_m<n$ and $t=s-1$ if $\gamma_m=n$. Then, the
$\k$-algebra ${\mathcal O}_q(G_{m,n}(\k))_\gamma$ is AS-Gorenstein if and only if (with
the above notation) $|\chi_{i-1}|=|\beta_i|$ for $1 \le i \le t$.
\end{subtheorem}

\proof As mentioned above, ${\mathcal O}_q(G_{m,n}(\k))_\gamma$ is a quantum graded
algebra with a straightening law on the poset
$\Pi_{m,n}\setminus\Pi_{m,n}^\gamma$. It follows from [LenRi(2); Remark 2.1.4]
that it has enough normal elements. 
 In addition,
${\mathcal O}_q(G_{m,n}(\k))_\gamma$ is a AS-Cohen-Macaulay domain, by Proposition
\ref{q-schubert-normal-domain} and Theorem \ref{q-sch-grass-CM}.
On the other hand,
${\mathcal O}_q(G_{m,n}(\k))_\gamma$ has a vector
space basis consisting of standard monomials on
$\Pi_{m,n}\setminus\Pi_{m,n}^\gamma$, since it is a 
quantum graded algebra with straightening law on this poset. 
Clearly, this implies that the Hilbert
series of ${\mathcal O}_q(G_{m,n}(\k))_\gamma$ is independent of the particular value of
$q \in\k^\ast$. Hence, using [JZ; Theorem 6.2], ${\mathcal O}_q(G_{m,n}(\k))_\gamma$ is
AS-Gorenstein if and only if ${\mathcal O}_1(G_{m,n}(\k))_\gamma$ is AS-Gorenstein; that
is, if and only if ${\mathcal O}_1(G_{m,n}(\k))_\gamma$ 
is Gorenstein (by the discussion above). The result now follows from 
Corollary 8.13 of [BV] . \qed

\section{Application to the normality of quantum determinantal rings.}
\label{maxor-qdetrings}

The aim of this section is to apply the previous results on quantum Schubert
varieties in order to show that quantum determinantal rings are normal
domains. Recall that a {\it quantum determinantal ring} is a factor
${\mathcal O}_q(M_{m,n}(\k))/\I_t$, where $1 \le t \le m$ and $\I_t$ is the ideal of
${\mathcal O}_q(M_{m,n}(\k))$ generated by the $t \times t$ quantum minors. \\

Let us start by defining a larger class of quotients of ${\mathcal
O}_q(M_{m,n}(\k))$. This class is obtained from ${\mathcal O}_q(M_{m,n}(\k))$
in the same way as quantum Schubert varieties are obtained from ${\mathcal
O}_q(G_{m,n}(\k))$.

\begin{definition} -- \label{def-gen-q-det-rings}
Let $\delta\in\Delta_{m,n}$ and put
$\Delta_{m,n}^\delta=\{\alpha\in\Delta_{m,n} \tq \alpha\not\ge_\st\delta\}$.
\[
{\mathcal O}_q(M_{m,n}(\k))_\delta 
= {\mathcal O}_q(M_{m,n}(\k))/\langle \Delta_{m,n}^\delta \rangle .
\]
\end{definition}

\begin{remark} -- \label{basic-rem-qdetgen} \\ \rm
(i) Let $\delta \in \Delta_{m,n}$. It is clear that $\Delta_{m,n}^\delta$ is a
$\Delta_{m,n}$-ideal. Hence, as mentioned in Subsection
\ref{q-gr-ASL-reminder}, the generalised quantum determinantal ring
${\mathcal O}_q(M_{m,n}(\k))_\delta$ inherits from ${\mathcal O}_q(M_{m,n}(\k))$ the structure of
a quantum graded A.S.L. on the poset
$\Delta_{m,n}\setminus\Delta_{m,n}^\delta$. Notice also that
$\Delta_{m,n}\setminus\Delta_{m,n}^\delta$ has a unique minimal element.\\
(ii) Let $\delta \in \Delta_{m,n}$. It follows from point (i) above that
${\mathcal O}_q(M_{m,n}(\k))_\delta$ has a standard monomial basis inhereted from the
corresponding basis of ${\mathcal O}_q(M_{m,n}(\k))$. The elements of this $\k$-basis
are the canonical images in ${\mathcal O}_q(M_{m,n}(\k))_\delta$ of the standard
monomials of ${\mathcal O}_q(M_{m,n}(\k))$ which either are equal to 
$1$ or are of the form
$[I_1|J_1]\dots[I_t|J_t]$ with $\delta \le_\st [I_1|J_1]$. \\
(iii) Let $\delta=(\{i_1 < \dots < i_r\},\{j_1 < \dots < j_r\})$, for some
integer $r$ such that $1 \le r \le m$. Then, 
\[ 
\GKdim {\mathcal O}_q(M_{m,n}(\k))_\delta = (m+n)r-\sum_{s=1}^r (i_s+j_s)+r, 
\] 
by Proposition
\ref{GKdim-q-gr-ASL} and [BV; 5.12].\\
(iv) As noticed
in [LenRi(2); \S 3.5], quantum determinantal rings are special cases of
generalised quantum determinantal rings, hence justifying the vocabulary in
Definition \ref{def-gen-q-det-rings}. More precisely, the paragraph before
[LenRi(2); Corollary 3.5.4] shows that, for $1 < t \le m$, the quantum 
determinantal ring 
${\mathcal O}_q(M_{m,n}(\k))/\I_t$ is equal to 
the generalised quantum determinantal
ring ${\mathcal O}_q(M_{m,n}(\k))_\delta$ where
$\delta=(\{1,\dots,t-1\},\{1,\dots,t-1\})$. \\
\end{remark}

The normality of quantum determinantal rings will be establised by applying
Proposition \ref{CS-normality-q-gr-ASL}. Hence, we first need to prove that
generalised quantum determinantal rings are integral domains. This is what we
do now. Fix an element $\delta \in \Delta_{m,n}$ and denote by $\gamma$ the
image of $\delta$ under the map $\delta_{m,n}$ introduced in Section
\ref{sec-basic-definitions}. Recall from Section 
 \ref{sec-basic-definitions} the dehomogenisation map
\[
\begin{array}{ccrcl}
D_{m,n} & : & {\mathcal O}_q(M_{m,n}(\k))[y,y^{-1};\phi] & \longrightarrow &
{\mathcal O}_q(G_{m,m+n}(\k))[[M]^{-1}] 
\end{array} 
\]
which sends $[I|J]$ to $[K_{(I,J)}][M]^{-1}$ and $y$ to $[M]$. 
It is clear that the ideal $\langle \Delta_{m,n}^\delta \rangle$ of
${\mathcal O}_q(M_{m,n}(\k))$ is stable under $\phi$. Hence, using the obvious
abuse of notation, there is a canonical isomorphism

\[
\begin{array}{rcl}
{\mathcal O}_q(M_{m,n}(\k))[y,y^{-1};\phi]/\langle \Delta_{m,n}^\delta
\rangle[y,y^{-1};\phi] & \cong & ({\mathcal O}_q(M_{m,n}(\k))/\langle
\Delta_{m,n}^\delta \rangle)[y,y^{-1};\phi] \cr & = &
{\mathcal O}_q(M_{m,n}(\k))_\delta\,[y,y^{-1};\phi].
\end{array}
\] 
(Here, the notation $\langle \Delta_{m,n}^\delta \rangle[y,y^{-1};\phi]$ on
the left hand side stands for the two-sided ideal of
${\mathcal O}_q(M_{m,n}(\k))[y,y^{-1};\phi]$ generated by $\Delta_{m,n}^\delta$.) 
Now, observe that the ideal $\langle \Pi_{m,n+m}^\gamma \rangle$ of
${\mathcal O}_q(G_{m,m+n}(\k))$ is completely prime, by Corollary
\ref{q-schubert-normal-domain}, and does not intersect the set of powers of
$[M]$, as one may easily prove using the standard basis of
${\mathcal O}_q(G_{m,m+n}(\k))$. Hence, there is a canonical isomorphism
\[
\begin{array}{rcl}
{\mathcal O}_q(G_{m,m+n}(\k))[[M]^{-1}]/\langle \Pi_{m,n}^\gamma \rangle[[M]^{-1}] &
\cong & ({\mathcal O}_q(G_{m,m+n}(\k))/\langle \Pi_{m,n}^\gamma
\rangle)[\overline{[M]}^{-1}] \cr & = &
{\mathcal O}_q(G_{m,m+n}(\k))_\gamma\,[\overline{[M]}^{-1}].
\end{array}
\]

(Recall that, if $x \in {\mathcal O}_q(G_{m,m+n}(\k))$, then $\overline{x}$ stands for
its canonical image in ${\mathcal O}_q(G_{m,m+n}(\k))_\gamma$. The same convention
applies to ${\mathcal O}_q(M_{m,n}(\k))_\delta$.)

Clearly, $D_{m,n}(\langle \Delta_{m,n}^\delta \rangle[y,y^{-1};\phi])=\langle
\Pi_{m,n}^\gamma \rangle)[[M]^{-1}]$; and so 
it follows that $D_{m,n}$ induces an
isomorphism
\[
\begin{array}{ccrcl}
D_{m,n}^\delta & : & {\mathcal O}_q(M_{m,n}(\k))_\delta[y,y^{-1};\phi] & \longrightarrow
& {\mathcal O}_q(G_{m,m+n}(\k))_\gamma\,[\overline{[M]}^{-1}] 
\end{array}.
\]

that sends $\overline{[I|J]}$ to 
$\overline{[K_{(I,J)}]}\,\overline{[M]}^{-1}$ 
and $y$ to $\overline{[M]}$. 

From this we deduce the following result.

\begin{proposition} -- \label{int-of-qdetgen}
Let $\delta \in \Delta_{m,n}$. Then the ring ${\mathcal O}_q(M_{m,n}(\k))_\delta$ is an
integral domain.
\end{proposition}

\proof Let $\gamma$ be the image of $\delta$ under the map $\delta_{m,n}$. By
the above discussion, we have an isomorphism
\[
D_{m,n}^\delta \,:\, {\mathcal O}_q(M_{m,n}(\k))_\delta\,[y,y^{-1};\phi] 
\longrightarrow
{\mathcal O}_q(G_{m,m+n}(\k))_\gamma\,[\overline{[M]}^{-1}]\,.
\]
Hence, the result follows by Corollary \ref{q-schubert-normal-domain}. \qed\\

The normality of quantum determinantal rings is now easy to obtain. First, we
need a lemma.

\begin{lemma} -- \label{lemma-qdet-loc}
Let $2 \le t \le m$. Let $\delta=[\{1,\dots,t-1\}|\{1,\dots,t-1\}]\in
{\mathcal O}_q(M_{m,n}(\Bbbk))$ and let $\overline{\delta}$ be the canonical
image of $\delta$ in ${\mathcal O}_q(M_{m,n}(\Bbbk))/{\mathcal I}_t$. 
Let  $A_t$ be the
subalgebra of ${\mathcal O}_q(M_{m,n}(\Bbbk))$ generated by the elements
$X_{ij}$ such that either $i \le t-1$ or $j \le t-1$. Then $\delta$ is a
normal element in $A_t$ and
\[
({\mathcal O}_q(M_{m,n}(\k))/{\mathcal I}_t)[\overline{\delta}^{-1}] \cong
A_t\,[\delta^{-1}]\,.
\].
\end{lemma}

\proof For $2 \le t \le m$, there is an obvious algebra homomorphism
\[
A_t 
\longrightarrow {\mathcal O}_q(M_{m,n}(\k))
\longrightarrow {\mathcal O}_q(M_{m,n}(\k))/\I_t
\longrightarrow ({\mathcal O}_q(M_{m,n}(\k))/\I_t)[\overline{\delta}^{-1}].
\]
Here, the first map is the obvious injection, the second is the
canonical projection while the third is the natural injection of a ring
into a localisation with respect to the multiplicative set generated
by a regular normal element. In addition,
$\delta$ is a normal element of $A_t$, see [PW; Lemma 4.5.1], and it 
is sent to an invertible element by the above map. Hence, there is 
a natural algebra homomorphism 
\[
\phi_t \, : \, 
A_t[\delta^{-1}] \longrightarrow ({\mathcal O}_q(M_{m,n}(\Bbbk))/\I_t)[\overline{\delta}^{-1}].
\]
We will show that $\phi_t$ is an isomorphism. First, we show that $\phi_t$ is
surjective. Clearly, $\overline{\delta}^{-1}$ is in the image of
$\phi_t$. Let $1 \le i \le m$ and $1 \le j \le n$. It is clear that if either
$i < t$ or $j < t$, then $X_{ij}+\I_t$ is in the image of $\phi_t$. Now,
assume $i,j \ge t$. The subalgebra of ${\mathcal O}_q(M_{m,n}(\Bbbk))$
generated by the entries of the generic matrix corresponding to rows
$1,\dots,t-1,i$ and columns $1,\dots,t-1,j$ is isomorphic to ${\mathcal
O}_q(M_t(\Bbbk))$. We may then develop its quantum determinant with respect to
its last row. This relation, seen inside $({\mathcal
O}_q(M_{m,n}(\Bbbk))/\I_t)[\overline{\delta}^{-1}]$, shows that $X_{ij}+\I_t$
is in the image of $\phi_t$.

Next we show that $\phi_t$ is injective. First, it is easy to see that 
\[
\GKdim A_t=mn-(m-(t-1))(n-(t-1)).
\] 
In addition, since $\delta$ commutes up to scalar
with each of the canonical generators of $A_t$, it is easy to show that it is
a local normal element in the sense of [KL; \S 12.4]. It then follows by [KL;
Theorem 12.4.4.] that $\GKdim (A_t)_{\delta}=mn-(m-(t-1))(n-(t-1))$. On the
other hand, by Remark \ref{basic-rem-qdetgen}, $\GKdim {\mathcal
O}_q(M_{m,n}(\Bbbk))/\I_t=mn-(m-(t-1))(n-(t-1))$ from which it follows that
$\GKdim ({\mathcal O}_q(M_{m,n}(\Bbbk))/{\mathcal
I}_t)[\overline{\delta}^{-1}] \ge mn-(m-(t-1))(n-(t-1))$, see [KL; Lemma
3.1]. As $A_t\,[\delta^{-1}]$ is an integral domain, any non-zero element
of $A_t\,[\delta^{-1}]$ is regular. Hence, in view of the above estimates of
Gelfand-Kirillov dimension, by [KL; Prop. 3.15] we must have
$\ker\phi_t=(0)$.\qed

\begin{theorem} -- \label{q-det-rings-are-normal-domains}
Let $1 \le t \le m$. Then ${\mathcal O}_q(M_{m,n}(\Bbbk))/\I_t$ is a normal domain.
\end{theorem}

\proof The case $t=1$ is trivial. Let $2 \le t \le m$. Recall that
${\mathcal O}_q(M_{m,n}(\k))/\I_t={\mathcal O}_q(M_{m,n}(\k))_\delta$, where
$\delta=(\{1,\dots,t-1\},\{1,\dots,t-1\})$. Hence, ${\mathcal O}_q(M_{m,n}(\k))/\I_t$ is
a domain by Proposition \ref{int-of-qdetgen}. In addition, we are in the
context of Section \ref{par-normality} since ${\mathcal O}_q(M_{m,n}(\k))/\I_t$ is a
quantum graded A.S.L. whose underlying poset has a single minimal element,
namely $\overline{\delta}=\delta + \I_t$. In addition, Proposition
\ref{int-of-qdetgen} shows that, for any upper neighbour $\gamma$ of
$\overline{\delta}$, the ring $({\mathcal O}_q(M_{m,n}(\k))/\I_t)/I_\gamma$ is a domain
(here, we are using the notation of Subsection \ref{par-normality}). By
Proposition \ref{CS-normality-q-gr-ASL}, it is enough to show that
$({\mathcal O}_q(M_{m,n}(\k))/\I_t)[\overline{\delta}^{-1}]$ is a normal domain.
However, using [MR; V. Proposition 2.5, IV. Proposition 2.1], this is a
consequence of Lemma \ref{lemma-qdet-loc}.\qed

\begin{remark} -- \rm In the light of Theorem 
\ref{q-det-rings-are-normal-domains}, 
one obvious question arises: are generalised quantum determinantal rings 
normal domains? Recall that this is true in the commutative case, see [BV;
Theorem 6.3]. Let $\delta\in\Delta_{m,n}$. The isomorphism $D_{m,n}^\delta$
together with Corollary \ref{q-schubert-normal-domain} shows that the skew
Laurent extension ${\mathcal O}_q(M_{m,n}(\k))_\delta\,[y,y^{-1};\phi]$ is normal.
However, we have not been able to deduce  that
${\mathcal O}_q(M_{m,n}(\k))_\delta$ is normal from this fact. 
Another approach to the normality of
${\mathcal O}_q(M_{m,n}(\k))_\delta$ (which is the approach we have used to derive the
normality of quantum determinantal rings) would be to apply Proposition
\ref{CS-normality-q-gr-ASL}. We would then need to show that
${\mathcal O}_q(M_{m,n}(\k))_\delta\,[\overline{\delta}^{-1}]$ is a normal domain. The
problem here is of technical nature: presumably,
${\mathcal O}_q(M_{m,n}(\k))_\delta\,[\overline{\delta}^{-1}]$ can be described as a
localisation of an iterated Ore extension (this would generalise Lemma
\ref{lemma-qdet-loc} to any generalised quantum determinantal rings). However,
getting such a description seems to be a rather tricky computation.
\end{remark}

\section*{References.}

{\bf [BV]} W. Bruns and U. Vetter. Determinantal rings. 
Lecture Notes in Mathematics, 1327. Springer-Verlag, Berlin, 1988. \newline
{\bf [JZ]} P. J\o rgensen and J.J. Zhang. {\it Gourmet's guide to Gorensteinness.}  
 Adv. Math.  151  (2000),  no. 2, 313--345. \newline
{\bf [GL]} N. Gonciulea and V. Lakshmibai. Flag Varieties. Actualit\'es Math\'ematiques. Hermann, Paris, 2001.\newline
{\bf [KLR]} A. Kelly, T.H. Lenagan and L. Rigal. {\em Ring theoretic properties of quantum grassmannians.} 
 J. Algebra Appl.  3  (2004),  no. 1, 9--30. \newline
{\bf [KL]} G. Krause and T.H. Lenagan. Growth of algebras and Gelfand-Kirillov dimension. 
Revised edition. Graduate Studies in Mathematics, 22. American Mathematical Society, Providence, RI, 2000. \newline
{\bf [KroLe]} D. Krob, and B. Leclerc. {\em Minor identities for quasi-determinants and quantum determinants.}  
Comm. Math. Phys.  169  (1995),  no. 1, 1--23. \newline
{\bf [LakRe]} V. Lakshmibai and N. Reshetikhin. {\em Quantum deformations of $SL_n/B$ and its Schubert varieties.} 
 Special functions (Okayama, 1990),  149--168, ICM-90 Satell. Conf. Proc., Springer, Tokyo, 1991. \newline
{\bf [LenRi(1)]} T.H. Lenagan and L. Rigal. {\em The maximal order property for quantum determinantal rings.}  
Proc. Edinb. Math. Soc. (2)  46  (2003),  no. 3, 513--529. \newline
{\bf [LenRi(2)]} T.H. Lenagan and L. Rigal. {\em Quantum graded algebras with a straightening law and the AS-Cohen-Macaulay property for quantum
determinantal rings and quantum grassmannians.}  J. Algebra  301  (2006),  no. 2, 670--702. \newline
{\bf [McCR]} J.C. McConnell and J.C. Robson. Noncommutative Noetherian rings. Revised edition. Graduate Studies in Mathematics, 30. 
American Mathematical Society, Providence, RI, 2001. \newline
{\bf [MR]} G. Maury and J. Raynaud. Ordres maximaux au sens de K. Asano. Lecture Notes in Mathematics, 808. Springer, Berlin, 1980. \newline
{\bf [PW]} B. Parshall and J. Wang. Quantum linear groups.  Mem. Amer. Math. Soc.  89  (1991),  no. 439. \newline

{\bf [R]} L. Rigal. {\em Normalit\'e de certains anneaux déterminantiels quantiques.} Proc. Edinburgh Math. Soc. (2)  42  (1999),  no. 3, 621--640. \newline

{\bf [Z]} J.J. Zhang.  {\em Connected graded Gorenstein algebras with enough normal elements.}  
J. Algebra  189  (1997),  no. 2, 390--405. 

\vskip 1cm
\noindent 
T. H. Lenagan: \\
Maxwell Institute for Mathematical Sciences\\
School of Mathematics, University of Edinburgh,\\
James Clerk Maxwell Building, King's Buildings, Mayfield Road,\\
Edinburgh EH9 3JZ, Scotland, UK\\~\\
E-mail: tom@maths.ed.ac.uk \\
\\
L. Rigal: \\
Universit\'e Jean-Monnet (Saint-\'Etienne), \\
Facult\'e des Sciences et Techniques, \\
D\'e\-par\-te\-ment de Math\'ematiques,\\
23 rue du Docteur Paul Michelon,\\
42023 Saint-\'Etienne C\'edex 2,\\France\\~\\
E-mail: Laurent.Rigal@univ-st-etienne.fr

\end{document}